\documentclass[10pt]{article}

\usepackage{xypic,amsfonts,amsmath,amssymb}
\usepackage[english]{babel}

\author{Luc Pirio}

\title{Study of a functional equation associated to the Kummer's
   equation of the trilogarithm. \\ Applications}  

\date{13/06/2002 }

\newtheorem{coro}{\textbf{corollary}}

\newtheorem{prop}{\textbf{Proposition}}
\def \abs  #1{{\left\vert #1\right\vert}}

\def \g {{\bf g}}

\def \h {{\bf h}}

\def \l  #1 {{ {\bf L}{\mbox{i}}_#1}}

\def \kk  {({\cal K}_3)}
\def \ccs { {{\mathbb C}^2 \setminus S}}
\def \ss { \underline{ {\cal S}}  }

\begin{document}

\maketitle

\begin{quote}
{\bf Abstract :} 
{\small In this paper we study a generalisation in 9 unknown
functions of a homogeneous version of the Kummer's equation for
$ \l{3} $. 
We give an explicit description of the space of local holomorphic solutions at 
 a generic point in $ \mathbb C ^2 $. Then we apply
this description firstly for obtaining new non linearisable maximal
rank webs (confirming some results annonced about one year ago by G. Robert
(\cite{hen1})). Secondly we show that under suitable conditions, the
trilogarithm is the only function which verifies the Kummer equation. }
\end{quote}

\section{Introduction}

\hspace{0.2cm} Since the works of Spence, Kummer and Abel, a big number of functional
equations satisfyed by low-order polylogarithms $ \l{n} $ are known ($n\leq 5$ ) ( see \cite{lewin}) .
For example, it was proved by Spence and independently by Kummer
that the trilogarithm $ \l{3} $ verifies the
fonctional equation 
\begin{align}
  2\l{3} (x) &  +  2 \l{3} (y) -   \l{3} (\frac{x}{y})   
+2 \l{3} ( \frac{1-x}{1-y} )  + 2 \l{3} (\frac{x(1-y)}{y(1-x)})    
-\l{3} (xy) \nonumber \\
 & +   2 \l{3} (-  \frac{x(1-y)}{(1-x)}          ) 
+ 2  \l{3} ( - \frac{(1-y)}{y(1-x)})   - \l{3} (\frac{x(1-y)^2}{y(1-x)^2}) 
 \qquad     \ ({\cal K}_{3})  \nonumber \\
  &=2\l{3}(1)-\log(y)^2\log( \frac{1-y}{1-x})
+\frac{\pi^2}{3}\log(y)+\frac{1}{3}\log(y)^3   \nonumber
\end{align}

for $x,y\in \mathbb R$ such that $ 0<x<y<1$ .\\
( from now on we note $E_3(x,y) $ the right member of $ \kk $ ). \\

Let us introduce the interior functions which appear in $ ({\cal
  K}_3)$: 
\begin{align}
U_1(x,y) & =x    &  U_2(x,y) & =y    &  U_3(x,y) & =\frac{x}{y}  \nonumber  \\
U_4(x,y) & =\frac{1-y}{1-x}   & U_5(x,y) & =
   \frac{x(1-y)}{y(1-x)}  &  U_6(x,y) & =xy   \nonumber  \\
U_7(x,y) & = \frac{x(1-y)}{x-1}   & U_8(x,y) & = \frac{1-y}{y(x-1)}  &  
   U_9(x,y) & =
   \frac{x(1-y)^2}{y(1-x)^2}         \nonumber   
\end{align}
They all are rational functions with real coefficients such that their
level curves in $\mathbb C ^2$ are lines, conics or cubics. \\

Associated to $\kk$ , we can consider the following homogeneous functional
equation in $9$ unknown functions: 
$$  F_1(U_1)+F_2(U_2)+F_3(U_3)+...+F_9(U_9)=0 \qquad \qquad ({\cal E})$$
In this paper our interest is in  the local solutions of this equation, a
local solution being a nine-uplet of function-germs 
satisfying the equation above.\\
 But we have to make it more precise.\\
Since we can consider equation $({\cal E})$ in any neighbourhood of any
$\omega\in \mathbb C ^2\setminus S'$ where $S'$ is the union of the
polar locus of the functions $U_i$ , we introduce ``the space of local
solutions in the class ${\cal F}$ of $({\cal E})$ at  $\omega$'': 
$$ \underline{{\cal
      S}{}_{\omega}^{\cal F}}= \left\{ {\bf
  { F}}=(F_i) \in \prod_{i=1}^9 \underline{{\cal
      F}_{\omega_i}} \; \: | \; \sum_1^9 F_i(U_i)=0 \: \text{ in } \:
  \underline{ {\cal F}_{\omega}} \: \right\} $$
(In this definition , $\omega_i$ denotes $U_i(\omega)$ for $i=1,..,9 $ ,
  ${\cal F}$ is any scheave of function-germs on $\mathbb K $ or $
\mathbb K ^2 $ ( with $\mathbb K =\mathbb R$ or $  \mathbb C$) , and $
\underline{{\cal F}_\theta}$ denotes the germ at $\theta\in \mathbb
K,\mathbb K ^2 $).\\

We will deal with at least measurable functions and to have non
pathologic situation, we will take $\omega$ generic in $\mathbb C ^2
\setminus S'$ more precisely such that 
the level curves of the $U_i$'s are not tangent in $\omega$.\\
Let be $S:=\bigl( \bigcup_{i<j} \left\{  \eta \in \mathbb C ^2 \: | \:
 dU_i \wedge dU_j(\eta)= 0 \: \right\} \bigr) \bigcup S'$ : from now
on, we will work with  the scheave of real measurable
function-germs noted $ {\cal M}$ and we will take $\omega\in \mathbb C^2 \setminus S $ .\\
Because the functions $U_i$ are rational functions, we can consider
that they are  defined on $ \mathbb C ^2 \setminus S$ and so the
equation $({\cal E})$ can be seen as a complex equation in the
complex field.\\
In the  part 2 of this note , we will consider the space of
holomorphic solutions of $ ({\cal E}) $ at $\omega \in \ccs $, i.e. the space 
$$    \underline{{\cal
      S}{}_{\omega}^{\cal O}}= \left\{ {\bf
  { F}}=(F_i) \in \prod_{i=1}^9 \underline{{\cal
      O}_{\omega_i}} \; \: | \; \sum_1^9 F_i(U_i)=0 \: \text{ in } \:
  \underline{ {\cal O}_{\omega}} \: \right\} $$
By a classical result of web geometry it comes that $   \underline{{\cal
      S}{}_{\omega}^{\cal O}}$ is a finite-dimensional $\mathbb
  C$-linear space and we have a bound $
  \text{dim}_{\scriptstyle{\mathbb C}} (  \underline{{\cal
      S}{}_{\omega}^{\cal O}} ) \leq 36 $ .\\
Next we will give a family $ {\bf \Gamma}$ of 36 linearily independant elements of 
$   \underline{{\cal S}{}_{\omega}^{\cal O}}$ . Thus  it comes that 
$ {\bf \Gamma}$ must be a basis and so it spans the whole space.\\

In part 3, we will apply the preceding results:\\
the fact that $   \underline{{\cal S}{}_{\omega}^{\cal O}}$ is of
maximal dimension 36 gives us some new examples of non linearisable
maximal rank planar webs. 
We will discuss this more precisely in part 3.1. ;
 in part 3.2. we will apply the explicit knowledge of $   \underline{{\cal
    S}{}_{\omega}^{\cal O}}$ given by  $ {\bf \Gamma}$ to the problem
of characterizing  $\l{3} $  by the equation $\kk$ of Kummer, what was
our initial goal.\\

{\bf Remark: } 
{\bf 1.} While I was working on the subject, I was told by G. Henkin 
that in a personnal communication to him
(\cite{hen1}), A. H\'enaut annonced that his colleague G. Robert 
had found that the Kummer's web is of maximal rank by constructing an
explicit basis of the space of abelian relations, what is equivalent
to part 2. of this paper. But no additional informations about this were
given until now.

{\bf 2.} This is a short version of a paper in preparation wich
 will display the results presented here in a more complete way as
well as some new results: we will show that any local
holomorphic solution of  $({\cal E})$ is ``a priori'' a global but
multiform solution, and use this to construct a method to
solve $({\cal E})$ by considering the monodromy  of those solutions.\\

{\bf Acknowledgments :} I would like to thank G. Henkin for
introducing me to this subject and useful discussions and 
 A. Bruter who helped me to put this paper in form.

\section{explicit resolution of $({\cal E})$  in the
  holomorphic case}

Let's take $\omega_0=(\frac{1}{3},\frac{1}{2})\in\mathbb R ^2 \setminus
S$. We will solve $({\cal E})$ at $\omega_0$ ( but we will find later
that the resolution we get gives a resolution at all $\omega'\in \mathbb C ^2
\setminus S$).\\
From now on, we note $\underline{{\cal S }}:=      \underline{{\cal
    S}{}_{\omega_0}^{\cal O}}$ .\\
The way that we will use to explicitly solve $({\cal E})$ is the
following:\\
we have this classical result of web geometry (see \cite{blabol}):
\begin{prop}
Let $N$ be a positive integer and $V_1,V_2,..,V_N $ be $N$  elements of
$\underline{{\cal O}_0}( \mathbb C ^2 , \mathbb C )$ such that we have
the generic condition $dV_i \wedge dV_j (0)\neq 0 $ ( $i<j$).\\
Then the space $ \Bigl\{ (G_i) \in \bigl( {\underline{{\cal O}_0}} \bigr)^N \: | \:
\sum^N G_i(V_i)=0  \: \Bigr\} $
has a finite dimension less than $\frac{N(N-1)}{2}$ .
\end{prop}

In the case of the equation $({\cal E})$ we succeed to
construct 36 linearily independant holomorphic soltions at $\omega_0$
. It shows that $\ss$ is of maximal possible dimension what is
exceptional (see part 3.1.) .\\

First we have to consider the constant solutions of $({\cal E})$. We
can easily construct them. They form a space of complex dimension
$8$. Let $\{ {\bf C}_i \}_{i=1,..,8}$ be a basis of it.\\

We have the 28 following 9-uplets of holomorphic germs. Verifiying
that they are 28 linearily independant elements of $\underline{{\cal
    S}}$ would be easy, but tedious, so we will skip this part.

\begin{align}
&{\bf F}_1 =  \biggl( {\bf L}og( \bullet ) ,- {\bf L}og( \bullet ),-
 {\bf L}og( \bullet ),0,0,0,0,0,0   \biggr) \nonumber \\ 
 & {\bf F}_2 =  \biggl( {\bf L}og(\frac{ \bullet }{1-\bullet}) ,0,{\bf
 L}og(\frac{1- \bullet }{\bullet} ) , - {\bf L}og(1- \bullet
 ),0,0,0,0,0   \biggr) \nonumber  \\
& {\bf F}_3 =  \biggl( {\bf L}og(1- \bullet ) ,- {\bf L}og(1- \bullet ),0,
 {\bf L}og( \bullet ),0,0,0,0,0  \biggr) \nonumber  \\
& {\bf F}_4 =  \biggl(0,0, {\bf L}og( \bullet ) , {\bf L}og( \bullet ),-
 {\bf L}og( \bullet ),0,0,0,0   \biggr) \nonumber   \\
& {\bf F}_5 =  \biggl(  {\bf L}og(1- \bullet ) ,  0,- {\bf L}og(1-
 \bullet ) ,0, {\bf L}og( 1-\bullet ),0,0,0,0   \biggr) \nonumber  \\
& {\bf F}_6 =  \biggl( {\bf L}og( \bullet ) , {\bf L}og( \bullet
 ),0,0,0,- {\bf L}og( \bullet ),0,0,0,  \biggr)  \nonumber \\
& {\bf F}_7 =  \biggl( {\bf L}og( \bullet ) ,0,0, {\bf L}og(
\bullet),0,0,-{\bf L}og( \bullet ) + i \pi ,0,0   \biggr) \nonumber     \\  
& {\bf F}_8 =  \biggl( \frac{1}{ \bullet}, 0 ,0,0,\frac{1}{
 \bullet},0,\frac{1}{ \bullet} -1, 0,0 \biggr)  \nonumber \\
& {\bf F}_9 =  \biggl( {\bf L}og(1- \bullet ) ,0,0,0,0 , - {\bf
  L}og(1- \bullet ), {\bf L}og(1- \bullet ),0,0  \biggr) \nonumber \\
& {\bf F}_{10} =  \biggl( 0,\bullet ,0,\bullet ,0,0,\bullet-1, 0,0
 \biggr)  \nonumber      \\
& {\bf F}_{11} =  \biggl(0,0,0,0,0, {\bf L}og(\bullet ) , -{\bf
  L}og(\bullet ),   {\bf L}og(\bullet ),0  \biggr) \nonumber   \\
& {\bf F}_{12} =  \biggl(  0 , {\bf L}og(\bullet ) ,0,0,0,0,-   {\bf
  L}og(1- \bullet ), {\bf L}og(1- \bullet ),0   \biggr) \nonumber  \\
 & {\bf F}_{13} =  \biggl( 0,0,0,0,0,0,{\bf L}og( \bullet ) ,{\bf L}og( \bullet ),-
 {\bf L}og( \bullet )- 2 i {\pi}   \biggr)&&  \nonumber   \\
& {\bf F}_{14} =  \biggl( 0,0,0,0,{\bf L}og(1- \bullet ) ,0, {\bf L}og(1-
\bullet ),0,- {\bf L}og(1- \bullet ) \biggr) \nonumber \\
& {\bf F}_{15} =  \biggl( 0,\frac{1}{\bullet} ,0,0,\bullet , 0,0, \bullet-1 ,0
\biggr)  \nonumber \\
& {\bf F}_{16} =  \biggl(  \bullet, 0 ,0, \frac{1}{
 \bullet},0,0,0,\frac{1}{ \bullet} -1 ,0 \biggr)  \nonumber 
\end{align}
\begin{align}
& {\bf F}_{17} = \biggl( 0,0,{\bf a} ( \bullet) ,0,0, -{\bf a} ( \bullet)
,0,0,-  {\bf a} ( \bullet) \biggr) \nonumber \\
& {\bf F}_{18} =  \biggl(2 {\bf L}og^2( \bullet ) ,2 {\bf L}og^2( \bullet ),-
 {\bf L}og^2( \bullet ),0,0,-{\bf L}og^2( \bullet ) ,0,0,0   \biggr)
 \nonumber  \\
& {\bf F}_{19} =  \biggl(0,0,0,0,0, {\bf L}og^2( \bullet ) ,-2 {\bf
  L}og( \bullet )^2,-2 {\bf L}og( \bullet )^2,{\bf L}og^2( \bullet
)+4i\pi{\bf L}og(\bullet)-4\pi ^2  \biggr) \nonumber \\
& {\bf F}_{20} =  \biggl(0,0, {\bf L}og^2( \bullet ) ,-2 {\bf L}og^2( \bullet ),-2{\bf L}og^2( \bullet ),0,0,0,{\bf L}og^2( \bullet )  \biggr)
 \nonumber  \\
& {\bf F}_{21} =  \biggl( {\bf d} ( \bullet ) ,-   {\bf d}(\bullet ),-
{\bf d}(\bullet ), -   {\bf d}(\bullet ),   {\bf d} ( \bullet ),
 0,0,0,0\biggr) \nonumber  \\
& {\bf F}_{22} =  \biggl(  {\bf d} ( \bullet ) ,   {\bf d}( \bullet )-
 \frac{i \pi }{2} \:  {\bf L}og(\bullet ),0,0,0,-
  {\bf d}(\bullet ),  {\bf d}   ( \bullet ), -  {\bf d}   ( \bullet ), 0\biggr)
 \nonumber \\
& {\bf F}_{23} =  \biggl( \pi ^2   ,0,0, {\bf d}   ( \bullet
  )-\frac{i\pi}{2}{\bf L}og(\bullet), {\bf d}  ( \bullet ),0,  {\bf d}   (\bullet ), {\bf d}( \bullet )+\frac{i\pi}{2}{\bf L}og(\bullet)-i\pi{\bf L}og(1-\bullet) ,-  {\bf d}    ( \bullet )\biggr) \nonumber \\
& {\bf F}_{24} =  \biggl( \l{2} ( \bullet ) , \l{2} ( \bullet
),0,\frac{1}{2} {\bf L}og ^2 (  \bullet),0,- \l{2} ( \bullet ),  \l{2}
( \bullet ),    -\l{2} ( \bullet ) -\frac{1}{2} {\bf L}og ^2 (
\bullet)  + i \pi {\bf L}og (
\bullet ) , \frac{\pi ^2 }{3}   \biggr)  \nonumber \\
& {\bf F}_{25} =  \biggl(0,\frac{1}{2} {\bf L}og ^2 (  \bullet), 0 ,
\l{2} (\bullet),\l{2} (\bullet),0,\l{2} (\bullet),\l{2}
(\bullet),-\l{2} (\bullet)  \biggr)  \nonumber \\
& {\bf F}_{26} =  \biggl(  2 \l{2} ( \bullet ), 0 
 ,- \l{2} ( \bullet ), 0 , 2 \l{2} ( \bullet ), - \l{2} (
  \bullet ),2 \l{2} ( \bullet ), 0 ,- \l{2} ( \bullet )
 \biggr)  \nonumber \\
&{\bf F}_{27}= \biggl(  2 {\g} (\bullet),2{\g}(\bullet),
-{\g}(\bullet),2{\g}(\bullet),2{\g}(\bullet),
-{\g}(\bullet),2\tilde{\g}(\bullet)
,2\tilde{\g}( \bullet),-{\g}(\bullet) \biggr) \nonumber \\
  &{\bf F}_{28}= \biggl(  2\h(\bullet),2\h(\bullet)-\frac{2 \pi ^2}{3}{\bf
  L}og( \bullet) ,
-\h(\bullet),2\h(\bullet),2\h(\bullet),
-\h(\bullet),2\tilde{\h}(\bullet)
,2\tilde{\h}( \bullet),- \h(\bullet) \biggr) \nonumber 
\end{align}
with
\begin{align} 
\bullet \; & := {\bf I}d_{\mathbb C} \nonumber \\
{\bf a} ( \bullet) & := {\bf a}\mbox{rcth} \:  ( \sqrt{ \bullet }\: )
\nonumber \\
{\bf d}(\bullet) & := \l{2}(\bullet)+\frac{1}{2}{\bf L}og(\bullet){\bf
  L}og(1-\bullet) -\frac{\pi ^2}{6} \nonumber \\
{\g}(\bullet) &  :=\l{3} (\bullet ) -{\bf L}og(\bullet) \; \l{2} (
\bullet )-\frac{1}{3}{\bf L}og^2( \bullet) \; {\bf
  L}og(1-\bullet)-\frac{2}{9}\l{3} (1) \nonumber \\
\h(\bullet)  &   := 2 {\bf L}og(\bullet) \; \l{2} (\bullet) + {\bf
  L}og^2 (\bullet) {\bf L}og(1-\bullet) \nonumber \\
\tilde{\g} (\bullet)  &  := \g(\bullet )-\frac{i\pi}{3}  \l{2}
  (\bullet) 
+\frac{4i\pi}{3}    {\bf
  d}( \bullet) +\frac{\pi ^2}{3}{\bf
  L}og(1-\bullet) +\frac{2i\pi ^3}{9} \nonumber \\
\tilde{\h}(\bullet)  &  := \h(\bullet) +2i \pi \l{2}
  (\bullet) -4i \pi {\bf d}(\bullet) -{\pi}^2
  {\bf L}og(1-\bullet) -\frac{2i\pi ^3}{3}\nonumber 
\end{align}
where all those functions are considered  holomorphic functions on
the whole simply connected domain $ \mathbb C \setminus 
\bigl( \{ 0 \} \times i \mathbb R^- \cup   \{ 1 \} \times i 
\mathbb R^+ \bigr) $, 
  functions which correspond to their usual definition on $]0,1[$ .\\

If we consider the family  ${\bf \Gamma }=\left\{ {\bf C}_i, {\bf F}_j \: | \: 
    1\leq i \leq 8 \: ,  1\leq j \leq 28 \: \right\} $  we obtain a family of 36 linearily
  independent elements of $\underline{\cal S}$ .  From the propostion 1, 
we know that  $\text{dim}_{ {}_{\scriptstyle{\mathbb C}}} \bigl(
  \underline{\cal S}  \bigr) \leq \frac{9(9-1)}{2}=36$  , and as
    we have seen before, this implies that  ${\bf \Gamma }$ is a basis
    of  $\underline{\cal S}$ .\\

So we have 
$$  \underline{ {\cal S}{}_{\omega_0}^{\cal O}}( {\cal E})=
\text{Vect}_{ {}_{\scriptstyle{ \mathbb C}}} \Bigl\langle \:  {\bf
  \Gamma }
\: \Bigr\rangle $$

{\bf Remark:}  Let $\omega'\in \ccs$ be different of $\omega_0$ .
There is a path $\gamma $ in $\ccs$ connecting $\omega_0$ to $\omega'$.
If ${\bf F}=(F_1,..,F_9)\in \ss $ we can easily see that each 
$F_i$ admits an analytic continuation along the path
$\gamma_i:=U_i\circ \gamma$ because this is verified for any element
of the basis ${\bf \Gamma}$.
It gives a holomorphic germ at $ \omega_i ':=U_i(\omega')$ noted
$F_i ^{[\gamma_i]}$. Then by analytic continuation along $\gamma$ and by the
unicity principle we get 
$$    \sum_{i=1}^9 F_i ^{[\gamma_i]} ( U_i ) =0 \qquad \text{in
  } \quad \underline{{\cal O}_{\omega'}} $$
and so ${\bf F}^{[\gamma]}:=(F_i ^{[\gamma_i]}) \in \underline{ {\cal
    S}{}_{\omega'}^{\cal O}}$ . It's clear that the application 
$ {\bf F} \rightarrow {\bf F}^{[\gamma]} $ is a linear isomorphism
between $  \underline{\cal S}$ and  $  \underline{ {\cal
    S}{}_{\omega'}^{\cal O}}$. This way, we can explicitly solve the
equation $({\cal E})$ at any point of $\ccs$ .

\section{Applications}
  
The fact that the dimension of $\ss$ is maximal and the explicit
description of $\ss$  both allow us to obtain some new results in
two a priori distinct subjects : the theory of planar webs
and the theory of polylogarithms.\\ 
For an introduction to the web theory, we refer to the
basic book of Blaschke and Bol `` Geometrie der Gewebe''
\cite{blabol}  and to \cite{cherngriff1} , \cite{chern} or \cite{web} 
for a more modern point of wiew.
As for polylogarithms, we refer to the books  \cite{lewin} and  \cite{pol} and to the talk of  J. Osterl\'e at the s\'eminaire Bourbaki (see  \cite{ost}).

\subsection{applications to web theory}
We suppose that the basic notions of web geometry are known.\\
From Bol's counterexample we know that not all the webs of maximal rank are
 linearisable (and so algebraic) : his web
 noted ${\cal B}$ is the global singular 5-web on $\mathbb C \mathbb P
 ^2$,  the 5 foliations of which are given by the level curves of the 
functions $U_1,U_2,..$ and $U_5$. Let's call $S_{\cal B}$ the singular locus
of $\cal B$ : it is the union of the polar locus of the $U_i$'s
($ i \leq 5$) with the algebraic set $ \cup_{1\leq i \leq 5} \{
w\in \mathbb C ^2  | \: dU_i\wedge U_j(w)=0 \}$~.  \\

It is a sub-web of the global singular 9-web noted ${\cal K }$ ( for
``Kummer'' ) defined by the level curves of the functions $U_i$ for
$i=1,..,9$ . Its singular locus is $S_{\cal K}=S$ . \\

From the elements $ {\bf
  F}_1,  {\bf F}_2, {\bf F}_3,  {\bf F}_4, {\bf F}_5$ and 
$   {\bf F}_{21}$  of $  \underline{\cal S}$ we can construct a base
  of the space ${\cal A}(\cal B)$ of the abelian relations of $\cal B$
  at $\omega_0$ .\\
So we have $\text{dim}_{\scriptstyle{\mathbb C}} {\cal A}({\cal B}) =6 $
  and  the web $\cal B$      is of maximal rank 6, although it is not 
linearisable.
From its discovery by Bol in the 30's onwards, this was the single
  known counterexample 
  to the problem of linearisation of planar webs of maximal rank.\\

From the fact that $\text{dim}_{\scriptstyle{\mathbb C}}   \underline{\cal
  S}=36$, we easily get that $\text{dim}_{\scriptstyle{\mathbb C}}
  {\cal A}({\cal K}) =28$, and because ${\cal K}$ is not linearisable ,
  ${\cal K}$  is another example of this kind of 
  web called ``exceptional planar webs '' by S.S. Chern.\\

According to Chern and Griffiths (see \cite{cherngriff2} page 83),
classifiying the non linearisable maximal rank webs is the fundamental
problem in web geometry. This explain the importance of this new
example of exceptional web.\\

But the explicit knowledge of the basis ${\bf \Gamma}$ of  $\underline{\cal
  S} $ allows to study all the sub-webs of ${\cal K}$ . For any subset $J\subset \{ 1,..,9\}$ we note ${\cal T}_J$
  the $|J|$-subweb of ${\cal K}$ given by the level curves of the
  function $U_j$, with $j\in J$. \\
If $ j_1,..,j_p$ are p distinct  integers in $ \{1,..,9\}$ , then we
  note $ \widehat{j_1..j_p}:=\{1,..,9\} \setminus \{ j_1,..,j_p\}$ .\\

\begin{prop}{\bf :}  \newline
\begin{tabular}{ll} 
$ \bullet$ ${\cal K} $  is an exceptional 9-web & $ \bullet$  ${\cal T}_{\widehat{ 69}} $
 is an exceptional 7-web \\
$ \bullet$  ${\cal T}_{\widehat{ 679}} $  is an exceptional 6-web & 
$ \bullet$  ${\cal T}_{\widehat{ 248}} $  is an exceptional 6-web \\
 \end{tabular}\\
Thoses two exceptional 6-webs are not equivalent ( up to local
  diffeomorphism~).\\
 $ {} \: \:  \bullet \: {\cal T}_{\widehat{ 369}}$  is a maximal-rank hexagonal
 6-web
\end{prop}
{\bf remarks :} {\bf 1.} By Bol's theorem ( see \cite{blabol} ), the fact that ${\cal
  T}_{\widehat{ 369}}$ is hexagonal implies that it is linearisable in
  a web formed by 6 pencils of lines. So it is algebraic, and the
  associated algebraic curve is an union of 6 lines in~$\mathbb C
  \mathbb P ^2$.\\
{\bf 2.}  The sub-webs ${\cal T}_{\widehat{ 36}} $ and 
${\cal T}_{\widehat{ 39}} $ are exceptional too but equivalent to
  ${\cal T}_{\widehat{ 69}} $ . \\
{\bf 3.} The sub-webs    ${\cal T}_{\widehat{ 689}} $, 
${\cal T}_{\widehat{ 349}} $,
${\cal T}_{\widehat{ 236}} $,
${\cal T}_{\widehat{ 359}} $, and
${\cal T}_{\widehat{ 136}} $
 are exceptional too
  but equivalent to ${\cal T}_{\widehat{ 679}} $.\\
 {\bf 4.} The sub-webs   ${\cal T}_{\widehat{ 147 }} $ ,
${\cal T}_{\widehat{ 257}} $ , and 
${\cal T}_{\widehat{ 158 }} $ 
are exceptional too
  but equivalent to ${\cal T}_{\widehat{248}} $.\\
 {\bf 5.} The two 6-webs ${\cal T}_{\widehat{ 679}} $ and ${\cal
  T}_{\widehat{ 248}} $ are not equivalent because one can prove
  that ${\cal T}_{\widehat{ 679}} $ contains an exceptional 5-subweb 
(the Bol's web ${\cal B}$!), contrarily to ${\cal
  T}_{\widehat{ 248}} $~, 5-subwebs of which have rank 5 and so are not exceptional.

{\bf 5.} We have a beautiful functional equation for $ \l{2} $ associated to
  ${\cal T}_{\widehat{248}} $ which is given by the element ${\bf
  F}_{26}$ of ${\bf \Gamma }$ :  \\ 
$ 2\l{2} (x) -\l{2}(\frac{x}{y})+ 2\l{2} (\frac{x(1-y)}{y(1-x)})
-\l{2}({x}{y}) + 2\l{2}
(-\frac{x(1-y)}{1-x})-\l{2}(\frac{x(1-y)^2}{y(1-x)^2 }) =0$ \\
(I have not seen an equation of this form in the bibliography).

\subsection{application to the caracterisation of $\l{3} $  by the \\
  equation $\kk$}

Our objective here is to study the function which satisfies the equation
$\kk$~.
This kind of problem has been studied for a long time for the Cauchy
equation $({\cal C})$: we know that any non-constant measurable
local solution of $({\cal C})$ is constructed from the logarithm.\\
There is similar results for the dilogarithm  (see
\cite{kies} and \cite{bloch}).\\
In his paper \cite{gon}, A. Goncharov obtains some results of
the same kind for the trilogarithm: \\

He considers the real single-valued cousin of $ \l{3} $  introduced by
Ramakhrishnan and Zagier :
$$  {\cal L}_3(z):=
 \Re e\left( \l{3} (z)-\log|z| \l{2} (z)+\frac{1}{3}\log|z|^2 \l{1}
 (z) \right) $$ 
defined on the whole $\mathbb C \mathbb P ^1$ and extended 
to $ \mathbb R [ \mathbb C \mathbb P ^1] $ by 
 linearity . \\
When it is well defined, he considers the following element of 
$ \mathbb Q [ \mathbb C \mathbb P ^1] $ :
\begin{align}
 R_3( \alpha_1, \alpha_2,\alpha_3):= \sum_{i=1}^3
& \biggl(  \{ \alpha_{i+2} \alpha_i -\alpha{i}+1 \} +   
             \{ \frac{\alpha_{i+2} \alpha_i -\alpha_i+1}{\alpha_{i+2}
               \alpha_i} \} +\{ \alpha_{i+2} \} \nonumber \\
&  + \{ \frac{\alpha_{i+2} \alpha_{i+1} -\alpha_{i+2}+1}{
   (\alpha_{i+2} \alpha_i -\alpha_{i}+1)    \alpha_{i+1} } \}  -     \{ \frac{\alpha_{i+2} \alpha_i
   -\alpha_i+1}{\alpha_{i+2}} \}     - \{ 1 \}
  \nonumber \\
&    - \{ \frac{\alpha_{i+2} \alpha_{i+1} -\alpha_{i+1}+1}{ (\alpha_{i+2} \alpha_i
   -\alpha_i+1) \alpha_{i+1} \alpha_{i+2}  } \}  + \{ \frac{\alpha_{i+2} \alpha_{i+1} -\alpha_{i+1}+1)\alpha_i}{ \alpha_{i+2} \alpha_i
   -\alpha_i+1 } \}     \biggr)    \nonumber \\             
& +  \{-\alpha_1 \alpha_2   \alpha_3 \}     \nonumber 
\end{align}
for $  \alpha_1,  \alpha_2,  \alpha_3 \in  \mathbb C \mathbb P ^1$
. (The indices i are taken modulo 3 ).\\
Next he proves that we have the functional equation in 22 terms 
$$ (\star \star ) \qquad \qquad {\cal L}_3 ( R_3( a,b,c))=0  \quad
\quad  \qquad a,b,c \in \mathbb C $$  

Then he shows ( part (a) of Theorem 1.10 in  \cite{gon} ) that \\
 {\it ``the space of real continuous functions on
$\mathbb C \mathbb P ^1\setminus \{ 0,1, \infty \}$ that satisfy the
functional equation $(\star \star )$ is generated by the functions
${\cal L}_3(z) $ and $ D_2(z).\log(|z|) $'' }  where $D_2$ is the
Bloch-Wigner function attached to $ \l{2} $  defined   for 
$ z\neq 0,1, \infty$ 
 by       $ D_2(z):=Im \left( \l{2} (z)+\log(1-z).\log(\abs{z}) \right) $ .\\

He had remarked before that if we specialize this equation by setting
$ a=1, b=x, $  and  $ c=\frac{1-y}{1-x} $ , the equation $(\star
\star)$  simplifies and by using the inversion relation $ {\cal
  L}_3({x}^{-1})= {\cal L}_3(x) \: , \: x  \in  \mathbb C \mathbb P ^1
$, it gives us exactly a homogeneous version
(i.e. without the second member $E_3(x,y)$  ) of the equation~$  \kk $~.

This leads him to ask if this specialization characterizes 
the solutions of $ (\star \star) $ .\\

The explicit determination of a basis of $\ss$ done in part 2. allows
us to give a positive answer to this question, in term of the function
$ \l{3} $ .\\
Let's first begin by a result of regularity for the mesurable
solutions of $({\cal E})$: 
\begin{prop}
Let be $\omega\in \mathbb R ^2 \setminus S$ and $ {\bf F}=(F_1,..,F_9)\in
\underline{ {\cal S}{}_{\omega}^{\cal M}}( {\cal E})$ .
Then each $F_i$ \\ is in fact an analytic germ at $\omega_i$ . Its 
complexification gives a germ $  F_i^{\scriptstyle{c}}  \in \underline{
  {\cal O}_{\omega_i}}$ such that $ {\bf F}^{\scriptstyle{c}}
:=(F_1^{\scriptstyle{c}},..,F_9^{\scriptstyle{c}})$ is a holomorphic
solution of $({\cal E})$ at~$\omega$~.
\end{prop}
{\bf sketch of the proof:} Because the level curves of the $U_i$'s are
in generic position near $\omega$, it comes from the paper of
A. Jar\`ai that the $F_i$'s are continuous germs ( see Theorem 3.3. in 
\cite{jarai} ). By elementary tools of integration you get next that
they are smooth germs. Then, similarly as Joly and Rauch in
\cite{jolyrauch} , one formulates the equation $({\cal E})$ in a
differential form. By an argument of ellipticity and by using
Petrowsky's theorem (see \cite{petro}) , you finally get that the
$F_i$'s are analytic germs . \\

So we have two $\mathbb R$-linear morphisms: \\
the first is just the restriction to $ \mathbb R ^2 $ with taking real
part
$$\begin{array}{rrcl}
   {\bf \rho} : &   \underline{ {\cal S}{}_{\omega}^{\cal O}}( {\cal E})  &
{\longrightarrow} & 
\underline{ {\cal S}{}_{\omega}^{\cal M}}( {\cal E} )  \nonumber \\ 
 \quad &  {\bf G} & \longmapsto &   \Re e( {\bf G}_{ | \scriptstyle{\mathbb
    R ^2} } )      
\end{array}$$
and the second is  given by the proposition 3 \\
$$ \begin{array}{rrcl}
 {\bf \varrho} : & \underline{ {\cal S}{}_{\omega}^{\cal M}} ( {\cal E})  &
{\longrightarrow }& 
\underline{ {\cal S}{}_{\omega}^{\cal O}}( {\cal E} )  \nonumber \\ 
 \quad &  {\bf F} & \longmapsto &   {\bf F}
  ^{\scriptstyle{c}} 
\end{array} $$

It is clear that $   {\bf \varrho } \circ  {\bf \rho} = {\bf I}d_{
  \scriptstyle{\underline{ {\cal S}{}_{\omega}^{\cal M}}( {\cal E}) }}     $ 
and so the study of measurable solutions of $({\cal E})$ at $\omega$ 
amounts to the study of the holomorphic
  solutions done in part 2. \\

We have this real semi-local characterization of $ \l{3} $ by  the
equation $ \kk $ :
\begin{prop}
Let $\epsilon_0$ be a real such that $ \frac{ \sqrt{5}-1}{2}<
\epsilon_0<1$  and \\ 
$ F  : ] -\infty,1 \: [ \: \rightarrow
\mathbb R $ be a measurable function such that for $ 0<x<y<\epsilon_0$ we have 
\begin{align}    
2 &  F(U_1(x,y))+  2 F(U_2(x,y))- F(U_3(x,y))      \nonumber \\
\qquad \qquad & \qquad \qquad \qquad   +  2 F(U_4(x,y)) +  2 F(U_5(x,y)) - F(U_6(x,y)   
\nonumber \\
 &  \qquad \qquad \qquad \qquad    + 2 F(U_7(x,y))+ 2 F(U_8(x,y))- F(U_9(x,y))=E_3(x,y) \nonumber 
\end{align}
\begin{itemize}  \item If $F$ is continuous at 0 then there
  exists $ a\in \mathbb R$ such that $$F:=\l{3}+a \: ( {\cal
  L}_3-\frac{2}{9}\l{3}( 1) )$$
\item If $F$ is derivable at 0 then  $ F:=\l{3} $ .
\end{itemize}
\end{prop}
{\ bf proof:} with our results of part 2. and the preceding remerk, it
is just a tedious exercise of linear algebra
The following statement, precisely related to Goncharov's question, is
equivalent to the preceding: 
\begin{coro}
Let $\epsilon_0$ be a real such that $ \frac{ \sqrt{5}-1}{2}<
\epsilon_0<1$  and \\ 
${\cal G} :   ] -\infty,1 [ \:       \rightarrow \mathbb R$  be
a measurable function such that for $ 0<x<y<\epsilon_0$ we have 
\begin{align}
 2 \: {\cal G}\left(x\right)&+2 \: {\cal G}\left(y\right) - \: {\cal
  G}\left(\frac{x}{y}\right) +2 \: {\cal
  G}\left(\frac{1-y}{1-x}\right) +2 \: {\cal G}\left(\frac{x(1-y)}{y(1-x)}\right) 
- {\cal G}(xy)  \nonumber \\  
& +2 \: {\cal G}\left( \frac{x(1-y)}{x-1}      \right)
+2  \:  {\cal G}\left( \frac{y-1}{y(1-x)}\right) -{\cal
  G}\left(\frac{x(1-y)^2}{y(1-x)^2} \right)= 2 \l{3} (1)    \nonumber
\end{align}
 Then if we suppose ${\cal G}$ continuous at $0$,then there existe $\alpha\in \mathbb R$ such that $$ {\cal G}=\alpha
 \: {\cal L}_3+\frac{2}{9}(1-\alpha) \: \l{3}(1) $$
\end{coro}

\begin{flushleft} {\it Luc Pirio, \\
Equipe d'analyse complexe,\\
Institut de mathématiques de Jussieu, \\
luclechat@hotmail.com\\
pirio@math.jussieu.fr}
\end{flushleft}


\begin{thebibliography}{99}
\selectlanguage{english}
\bibitem[Bla-Bo]{blabol} W. Blaschke, G. Bol, {\em Geometrie der Gewebe},
    Springer, Berlin, 1938 
\bibitem[Blo]{bloch} S. Bloch , {\em Higher regulators, algebraic
    K-thory and zeta functions of elliptic curves}, CRM Monograph
  Series {\bf 11} , A.M.S. providence, 2000




\bibitem[Ch-Gr1]{cherngriff1} S.S. Chern and  P.A. Griffiths , {\em   Abel's Theorem and Webs},
  Jahresber. Deutsch. Math.-Verein. {\bf 80} (1978), p. 13-110


\bibitem[Ch-Gr2]{cherngriff2} S.S. Chern and  P.A. Griffiths , {\em
    Correction and Addenda to our Paper: Abel's Theorem and Webs},
  Jahresber. Deutsch. Math.-Verein. {\bf 83} (1981), p. 78-83

\bibitem[Che]{chern} S.S. Chern, {\em Web Geometry },
  Bull. Amer. Math. soc. {\bf 6} (1982), p. 1-8


\bibitem[Gon]{gon}  Goncharov, A. , {\em Geometry of configurations,
    Polylogarithms and Motivic Cohomology}, Advances in Maths. {\bf 114}
    (1995) , p. 197-318 
%



\bibitem[H\'e1]{hen1} A. H\'enaut,  personnal letter to G. Henkin , 
  15 november 2001 

\bibitem[Jar]{jarai} J\'arai, A. , {\em On regular solutions of
    functional equations} , Aequationes Math. {\bf 30}, (1986), p. 21-54

\bibitem[Jo-Ra]{jolyrauch} J.L. Joly, J. Rauch ,{\em Ondes oscillantes
    semi-lin\'eaires en une dimension }~ Journ\'ees `` Equations aux
    d\'eriv\'ees partielles''( Saint Jean de Monts, 1986), Ecole
    polytechnique, pailaiseau, 1986

\bibitem[Kie]{kies} H. Kieswetter {\em Eine Bemerkung \"uber partielle
    Differentiationen bei N.H. Abel }, Publ. Math. Debrecen {\bf 5}
    (1957), p. 265-268


\bibitem[Lew]{lewin} L. {Lewin},  {\em Polylogarithms and
    Associated Functions } ,
 Elsevier North-Holland,
 New-York, 1981

\bibitem[Ost]{ost} J. Oesterl{\'e}, {\em Polylogarithmes }, S{\'e}minaire
    BOURBAKI, Vol 1992/93, Ast\'erisque {\bf 216}, Exp. n° 762, p.49-67


\bibitem[Pet]{petro} Petrowsky, I. {\em Sur l'analycit\'e des solutions des
    syst\`emes d'\'equations diff\'erentielles} , Math. Sbornik, {\bf
    47} (1939), p. 3-70
\bibitem[Pol]{pol} (L. {Lewin} editor) {\em Structural Properties of
    Polylogarithms  } , Maths. Surveys and Monographs, Vol. 37,  1991.
\bibitem[Web]{web} ( J. Grifone , E. Salem  editors ) {\em Web theory
    and related topics }, World scientific, 2001





\end{thebibliography}
\end{document}